\documentclass[preprint,times,12pt]{elsarticle}
\usepackage{tikz}

\newtheorem{conjecture}{Conjecture}
\newtheorem{lemma}{Lemma}
\newtheorem{fact}{Fact}
\newtheorem{corollary}{Corollary}
\newtheorem{theorem}{Theorem}

\newcommand{\Proof}{\noindent \textbf{Proof. }}
\newcommand{\sval}[1]{\left\llbracket #1 \right\rrbracket}
\newcommand{\dmin}[1]{d_{min}(#1)}
\newcommand{\dmax}[1]{d_{max}(#1)}

\newcommand{\FigureRho}{
\begin{figure}
\begin{center}
\begin{tikzpicture}[scale = 1.4]
\draw (0,0) circle(0.4mm);
\node at (0,0.4) {\( \frac{1}{2} \)};

\draw[fill = red] (0.24,0) circle(0.4mm);
\draw[color = red] (0.28,0) -- (2.42,0);
\node at (1.34,0.4) {\( \left[\frac{1}{2}{+}2^{-n}, \frac{2}{3}{-}\frac{7}{3}2^{-n}\right]\)};
\node at (1.5,-0.15) {\(  \delta = 2^{-n+1}\)};
\draw[fill = red] (2.46,0) circle(0.4mm);

\draw[fill = black] (2.92,0) circle(0.4mm);
\node at (2.92,1) {\(  \frac{2}{3}{-}\frac{1}{3}2^{-n}\)};
\draw (2.92,0) -- (2.92,0.7);

\draw (3,0) circle(0.4mm);
\node at (3,0.4) {\( \frac{2}{3} \)};

\draw[fill = green] (3.40,0) circle(0.4mm);
\draw[color = green] (3.44,0) -- (8.72,0);
\node at (6.08,0.4) {\( \left[\frac{2}{3}{+}\frac{5 }{3}2^{-n}, 1{-}2^{-n}\right]\)};
\node at (6.08,-0.15) {\( \delta = 2^{-n+1}\)};
\draw[fill = green] (8.76,0) circle(0.4mm);

\draw (9,0) circle(0.4mm);
\node at (9,0.4) {\(  1 \)};

\draw[fill = black] (0,-2) circle(0.4mm);
\node at (0,-2.4) {\( \frac{1}{2} \)};

\draw[fill = green] (0.36,-2) circle(0.4mm);
\draw[-stealth,color = green] (3.38,-0.02) -- (0.38,-1.98);
\draw[color = green] (0.40,-2) -- (4.28,-2);
\node at (2.34,-2.4) {\( \left[\frac{1}{2}{+}3\cdot 2^{-n-1}, \frac{3}{2^2}{-}2^{-n-1}\right]\)};
\node at (2.34,-1.85) {\( \delta = 3\cdot 2^{-n-1}\)};
\draw[fill = green] (4.38,-2) circle(0.4mm);
\draw[-stealth,color = green] (8.74,-0.02) -- (4.40,-1.98);

\draw (4.5,-2) circle(0.4mm);
\node at (4.5,-2.4) {\( \frac{3}{2^2}  \)};

\draw[fill = red] (4.98,-2) circle(0.4mm);
\draw[-stealth,color = red] (0.26,-0.02) -- (4.96,-1.98);
\draw[color = red] (5.02,-2) -- (8.24,-2);
\node at (6.51,-2.4) {\( \left[\frac{3}{2^2}{+}2^{-n+1}, 1-3{\cdot}2^{-n} \right]\)};
\node at (6.51,-1.85) {\( \delta = 3\cdot 2^{-n}\)};
\draw[fill = red] (8.28,-2) circle(0.4mm);
\draw[-stealth,color = red] (2.46,-0.02) -- (8.26,-1.98);

\draw (9,-2) circle(0.4mm);
\node at (9,-2.4) {\(  1 \)};

\draw[-stealth] (2.90,-0.02) -- (0.02,-1.97);

\node at (8.65,-1.2) {\small odd \( n \)};

\end{tikzpicture}
\vspace{1.5cm}

\begin{tikzpicture}[scale = 1.4]
\draw (0,0) circle(0.4mm);
\node at (0,0.4) {\( \frac{1}{2} \)};

\draw[fill = red] (0.24,0) circle(0.4mm);
\draw[color = red] (0.28,0) -- (2.56,0);
\node at (1.34,0.4) {\( \left[\frac{1}{2}{+}2^{-n}, \frac{2}{3}{-}\frac{5}{3}2^{-n}\right]\)};
\node at (1.5,-0.15) {\(  \delta = 2^{-n+1}\)};
\draw[fill = red] (2.60,0) circle(0.4mm);

\draw (3,0) circle(0.4mm);
\node at (3,0.4) {\( \frac{2}{3} \)};

\draw[fill = green] (3.08,0) circle(0.4mm);
\draw[color = green] (3.12,0) -- (8.72,0);
\node at (6.08,0.4) {\( \left[\frac{2}{3}{+}\frac{1}{3}2^{-n}, 1{-}2^{-n}\right]\)};
\node at (6.08,-0.15) {\( \delta = 2^{-n+1}\)};
\draw[fill = green] (8.76,0) circle(0.4mm);

\draw (9,0) circle(0.4mm);
\node at (9,0.4) {\( 1 \)};

\draw (0,-2) circle(0.4mm);
\node at (0,-2.4) {\( \frac{1}{2} \)};

\draw[fill = green] (0.12,-2) circle(0.4mm);
\draw[-stealth,color = green] (3.06,-0.02) -- (0.14,-1.98);
\draw[color = green] (0.16,-2) -- (4.34,-2);
\node at (2.34,-2.4) {\( \left[\frac{1}{2}{+}2^{-n-1}, \frac{3}{2^2}{-}2^{-n-1}\right]\)};
\node at (2.34,-1.85) {\( \delta = 3\cdot 2^{-n-1}\)};
\draw[fill = green] (4.38,-2) circle(0.4mm);
\draw[-stealth,color = green] (8.74,-0.02) -- (4.40,-1.98);

\draw (4.5,-2) circle(0.4mm);
\node at (4.5,-2.4) {\( \frac{3}{2^2} \)};

\draw[fill = red] (4.98,-2) circle(0.4mm);
\draw[-stealth,color = red] (0.26,-0.02) -- (4.96,-1.98);
\draw[color = red] (5.02,-2) -- (8.48,-2);
\node at (6.51,-2.4) {\( \left[\frac{3}{2^2}{+}2^{-n+1}, 1{-}2^{-n+1} \right]\)};
\node at (6.51,-1.85) {\( \delta = 3\cdot 2^{-n}\)};
\draw[fill = red] (8.52,-2) circle(0.4mm);
\draw[-stealth,color = red] (2.62,-0.02) -- (8.50,-1.98);

\draw (9,-2) circle(0.4mm);
\node at (9,-2.4) {\( 1 \)};

\node at (8.5,-1.2) {\small even \( n\)};

\end{tikzpicture}
\end{center}
\caption{\label{FigureRho} The function \( \rho(y) \) for odd \( n \) and for even \( n \). The \( y \) are on the top line and the \( \rho(y) \) are on the bottom line; the \( \rho(y) \) with \( h = 1 \) are red, those with \( h = 2 \) are green and those such that \( \rho(y) = 1/2 \) are black. The first value, the last value and the step \( \delta \) are annotated on each homogeneous interval of \( y \) and \( \rho(y) \).
}
\end{figure}
}

\newcommand{\CirclePoints}[3]{
\draw (0,0) circle(#2);
\foreach \ind in {0,1,...,#3}{
  \pgfmathparse{(mod((\ind)*#1,1)) *360}
  \node at (\pgfmathresult:#2+0.3) {\ind};
  \draw[fill = black] (\pgfmathresult:#2) circle(0.4mm);
  }
}

\newcommand{\XkCircle}{
\begin{figure}
\begin{center}
\begin{tikzpicture}
\draw[fill = black] (0,0) circle(0.4mm);
\node at (0,-1) {\(z = \log_2 3 \)};
\CirclePoints{log2(3)}{2.5}{15}
\end{tikzpicture} 
\end{center}
\caption{\label{XkCircle} Representation of the set \( X_{16} \), for \(z = \log_2 3 \), on a circle of lenght \( 1 \). Points \( x \in X_k \) are labeled by \( i \) such that \( x = (iz) \bmod 1 \).}
\end{figure}
}

\newcommand{\BoundDmax}{
\begin{figure}
\begin{center}
\begin{tikzpicture}
\draw[-stealth] (0,0)--(12,0);
\node at (12.2,0) {\( k \)};
\foreach \t in {0,0.5,...,11.5}{
  \draw[dotted] (\t,0)--(\t,3.5);
  }
\node at (0,-0.3) {\(k_8 \)};
\node at (11.5,-0.3) {\(k_9 \)};

\draw[fill=black] (0.0,2.342)circle(0.04)--(0.48,2.342);
\draw[fill=black] (0.5,2.242)circle(0.04)--(0.98,2.242);
\draw[fill=black] (1.0,2.142)circle(0.04)--(1.48,2.142);
\draw[fill=black] (1.5,2.042)circle(0.04)--(1.98,2.042);
\draw[fill=black] (2.0,1.942)circle(0.04)--(2.48,1.942);
\draw[fill=black] (2.5,1.842)circle(0.04)--(2.98,1.842);
\draw[fill=black] (3.0,1.742)circle(0.04)--(3.48,1.742);
\draw[fill=black] (3.5,1.642)circle(0.04)--(3.98,1.642);
\draw[fill=black] (4.0,1.542)circle(0.04)--(4.48,1.542);
\draw[fill=black] (4.5,1.442)circle(0.04)--(4.98,1.442);
\draw[fill=black] (5.0,1.342)circle(0.04)--(5.48,1.342);
\draw[fill=black] (5.5,1.242)circle(0.04)--(5.98,1.242);
\draw[fill=black] (6.0,1.142)circle(0.04)--(6.48,1.142);
\draw[fill=black] (6.5,1.042)circle(0.04)--(6.98,1.042);
\draw[fill=black] (7.0,0.942)circle(0.04)--(7.48,0.942);
\draw[fill=black] (7.5,0.842)circle(0.04)--(7.98,0.842);
\draw[fill=black] (8.0,0.742)circle(0.04)--(8.48,0.742);
\draw[fill=black] (8.5,0.642)circle(0.04)--(8.98,0.642);
\draw[fill=black] (9.0,0.542)circle(0.04)--(9.48,0.542);
\draw[fill=black] (9.5,0.442)circle(0.04)--(9.98,0.442);
\draw[fill=black] (10.0,0.342)circle(0.04)--(10.48,0.342);
\draw[fill=black] (10.5,0.242)circle(0.04)--(10.98,0.242);
\draw[fill=black] (11.0,0.142)circle(0.04)--(11.48,0.142);
\draw[fill=black] (11.5,0.1)circle(0.04)--(11.7,0.1);

\draw[color=red,fill=red] (0.0,3.270)circle(0.02)--(0.48,3.270);
\draw[color=red,fill=red] (0.5,3.134)circle(0.02)--(0.98,3.134);
\draw[color=red,fill=red] (1.0,2.998)circle(0.02)--(1.48,2.998);
\draw[color=red,fill=red] (1.5,2.862)circle(0.02)--(1.98,2.862);
\draw[color=red,fill=red] (2.0,2.725)circle(0.02)--(2.48,2.725);
\draw[color=red,fill=red] (2.5,2.589)circle(0.02)--(2.98,2.589);
\draw[color=red,fill=red] (3.0,2.453)circle(0.02)--(3.48,2.453);
\draw[color=red,fill=red] (3.5,2.317)circle(0.02)--(3.98,2.317);
\draw[color=red,fill=red] (4.0,2.180)circle(0.02)--(4.48,2.180);
\draw[color=red,fill=red] (4.5,2.044)circle(0.02)--(4.98,2.044);
\draw[color=red,fill=red] (5.0,1.908)circle(0.02)--(5.48,1.908);
\draw[color=red,fill=red] (5.5,1.772)circle(0.02)--(5.98,1.772);
\draw[color=red,fill=red] (6.0,1.635)circle(0.02)--(6.48,1.635);
\draw[color=red,fill=red] (6.5,1.499)circle(0.02)--(6.98,1.499);
\draw[color=red,fill=red] (7.0,1.363)circle(0.02)--(7.48,1.363);
\draw[color=red,fill=red] (7.5,1.226)circle(0.02)--(7.98,1.226);
\draw[color=red,fill=red] (8.0,1.090)circle(0.02)--(8.48,1.090);
\draw[color=red,fill=red] (8.5,0.954)circle(0.02)--(8.98,0.954);
\draw[color=red,fill=red] (9.0,0.818)circle(0.02)--(9.48,0.818);
\draw[color=red,fill=red] (9.5,0.681)circle(0.02)--(9.98,0.681);
\draw[color=red,fill=red] (10.0,0.545)circle(0.02)--(10.48,0.545);
\draw[color=red,fill=red] (10.5,0.409)circle(0.02)--(10.98,0.409);
\draw[color=red,fill=red] (11.0,0.272)circle(0.02)--(11.48,0.272);
\draw[color=red,fill=red] (11.5,0.195)circle(0.02)--(11.7,0.195);

\draw[color=red,fill=red] (0.0,1.635)circle(0.02)--(0.48,1.635);
\draw[color=red,fill=red] (0.5,1.564)circle(0.02)--(0.98,1.564);
\draw[color=red,fill=red] (1.0,1.493)circle(0.02)--(1.48,1.493);
\draw[color=red,fill=red] (1.5,1.422)circle(0.02)--(1.98,1.422);
\draw[color=red,fill=red] (2.0,1.351)circle(0.02)--(2.48,1.351);
\draw[color=red,fill=red] (2.5,1.280)circle(0.02)--(2.98,1.280);
\draw[color=red,fill=red] (3.0,1.209)circle(0.02)--(3.48,1.209);
\draw[color=red,fill=red] (3.5,1.138)circle(0.02)--(3.98,1.138);
\draw[color=red,fill=red] (4.0,1.066)circle(0.02)--(4.48,1.066);
\draw[color=red,fill=red] (4.5,0.995)circle(0.02)--(4.98,0.995);
\draw[color=red,fill=red] (5.0,0.924)circle(0.02)--(5.48,0.924);
\draw[color=red,fill=red] (5.5,0.853)circle(0.02)--(5.98,0.853);
\draw[color=red,fill=red] (6.0,0.782)circle(0.02)--(6.48,0.782);
\draw[color=red,fill=red] (6.5,0.711)circle(0.02)--(6.98,0.711);
\draw[color=red,fill=red] (7.0,0.640)circle(0.02)--(7.48,0.640);
\draw[color=red,fill=red] (7.5,0.569)circle(0.02)--(7.98,0.569);
\draw[color=red,fill=red] (8.0,0.498)circle(0.02)--(8.48,0.498);
\draw[color=red,fill=red] (8.5,0.427)circle(0.02)--(8.98,0.427);
\draw[color=red,fill=red] (9.0,0.355)circle(0.02)--(9.48,0.355);
\draw[color=red,fill=red] (9.5,0.284)circle(0.02)--(9.98,0.284);
\draw[color=red,fill=red] (10.0,0.214)circle(0.02)--(10.48,0.214);
\draw[color=red,fill=red] (10.5,0.142)circle(0.02)--(10.98,0.142);
\draw[color=red,fill=red] (11.0,0.071)circle(0.02)--(11.48,0.071);
\draw[color=red,fill=red] (11.5,0.090)circle(0.02)--(11.7,0.090) ;
\end{tikzpicture} 
\end{center}
\caption{\label{BoundDmax} The value of \( \dmax{k} \) (black) in the range \( k_8 \leq k < k_9 \) for \( z = \log_2 3 \) and the bounds (red) of Lemma \ref{dmaxkBound}. The range \( k_8 \leq k < k_9 \) is partitioned into \( q_{8} = 23 \) segments of \( \ell_{8} = 655 \) elements each.}
\end{figure}
}

\newcommand{\BoundkDmax}{
\begin{figure}
\begin{center}
\begin{tikzpicture}
\draw[-stealth] (0,0)--(12,0);
\node at (12.2,0) {\( k \)};
\foreach \t in {0,0.5,...,11.5}{
  \draw[dotted] (\t,0)--(\t,5);
  }
\node at (0,-0.3) {\(k_8 \)};
\node at (11.5,-0.3) {\(k_9 \)};

\draw[fill=black] (0.0,0.477)circle(0.05)--(0.48,0.804);
\draw[fill=black] (0.5,0.770)circle(0.05)--(0.98,1.082);
\draw[fill=black] (1.0,1.035)circle(0.05)--(1.48,1.333);
\draw[fill=black] (1.5,1.271)circle(0.05)--(1.98,1.556);
\draw[fill=black] (2.0,1.480)circle(0.05)--(2.48,1.751);
\draw[fill=black] (2.5,1.661)circle(0.05)--(2.98,1.918);
\draw[fill=black] (3.0,1.814)circle(0.05)--(3.48,2.057);
\draw[fill=black] (3.5,1.939)circle(0.05)--(3.98,2.168);
\draw[fill=black] (4.0,2.036)circle(0.05)--(4.48,2.251);
\draw[fill=black] (4.5,2.105)circle(0.05)--(4.98,2.306);
\draw[fill=black] (5.0,2.147)circle(0.05)--(5.48,2.334);
\draw[fill=black] (5.5,2.160)circle(0.05)--(5.98,2.333);
\draw[fill=black] (6.0,2.145)circle(0.05)--(6.48,2.304);
\draw[fill=black] (6.5,2.103)circle(0.05)--(6.98,2.248);
\draw[fill=black] (7.0,2.032)circle(0.05)--(7.48,2.164);
\draw[fill=black] (7.5,1.934)circle(0.05)--(7.98,2.051);
\draw[fill=black] (8.0,1.808)circle(0.05)--(8.48,1.911);
\draw[fill=black] (8.5,1.654)circle(0.05)--(8.98,1.743);
\draw[fill=black] (9.0,1.472)circle(0.05)--(9.48,1.547);
\draw[fill=black] (9.5,1.262)circle(0.05)--(9.98,1.323);
\draw[fill=black] (10.0,1.024)circle(0.05)--(10.48,1.071);
\draw[fill=black] (10.5,0.758)circle(0.05)--(10.98,0.791);
\draw[fill=black] (11.0,0.464)circle(0.05)--(11.48,0.484);

\draw[color=red,fill=red] (0.0,0.667)circle(0.05)--(0.48,1.123);
\draw[color=red,fill=red] (0.5,1.275)circle(0.05)--(0.98,1.711);
\draw[color=red,fill=red] (1.0,1.826)circle(0.05)--(1.48,2.242);
\draw[color=red,fill=red] (1.5,2.319)circle(0.05)--(1.98,2.715);
\draw[color=red,fill=red] (2.0,2.754)circle(0.05)--(2.48,3.130);
\draw[color=red,fill=red] (2.5,3.130)circle(0.05)--(2.98,3.487);
\draw[color=red,fill=red] (3.0,3.450)circle(0.05)--(3.48,3.786);
\draw[color=red,fill=red] (3.5,3.710)circle(0.05)--(3.98,4.027);
\draw[color=red,fill=red] (4.0,3.913)circle(0.05)--(4.48,4.210);
\draw[color=red,fill=red] (4.5,4.058)circle(0.05)--(4.98,4.335);
\draw[color=red,fill=red] (5.0,4.145)circle(0.05)--(5.48,4.403);
\draw[color=red,fill=red] (5.5,4.174)circle(0.05)--(5.98,4.412);
\draw[color=red,fill=red] (6.0,4.145)circle(0.05)--(6.48,4.363);
\draw[color=red,fill=red] (6.5,4.058)circle(0.05)--(6.98,4.256);
\draw[color=red,fill=red] (7.0,3.913)circle(0.05)--(7.48,4.091);
\draw[color=red,fill=red] (7.5,3.710)circle(0.05)--(7.98,3.869);
\draw[color=red,fill=red] (8.0,3.449)circle(0.05)--(8.48,3.558);
\draw[color=red,fill=red] (8.5,3.130)circle(0.05)--(8.98,3.249);
\draw[color=red,fill=red] (9.0,2.754)circle(0.05)--(9.48,2.853);
\draw[color=red,fill=red] (9.5,2.319)circle(0.05)--(9.98,2.398);
\draw[color=red,fill=red] (10.0,1.826)circle(0.05)--(10.48,1.886);
\draw[color=red,fill=red] (10.5,1.275)circle(0.05)--(10.98,1.315);
\draw[color=red,fill=red] (11.0,0.667)circle(0.05)--(11.48,0.729);

\draw[color=red] (0.02,4.717)--(11.48,4.701);

\end{tikzpicture} 
\end{center}
\caption{\label{BoundkDmax} The values of \( k\dmax{k} \) in the range \( k_8 \leq k < k_9 \) for \( z = \log_2 3 \). 
The range \( k_8 \leq k < k_9 \) is partitioned into \( q_{k_8} = 23 \) segments of \( \ell_{k_8} = 655 \) elements each. 
Red lines represent the upper bound \( 2\frac{q_{k_h}-t}{q_{k_h}}(1+t+\frac{j}{k_h}) \) from Lemma \ref{epsilonkpkt} and the global upper bound \( \frac{q_{k_h}+5+5/q_{k_h}}{2} \) from Lemma \ref{epsilonkpk}.}
\end{figure}
}

\title{Some contributions to Collatz conjecture.}

\author{Livio Colussi}
\address{Department of Mathematics\\
University of Padova\\
via Trieste, 63\\
35121 Padova (Italy)\\
e-mail: colussi@math.unipd.it
}

\begin{document}
\begin{abstract}
 The Collatz conjecture, also known as \( 3x+1 \) conjecture, can be stated in terms of the reduced Collatz function \( R(x) = (3x+1)/2^m \) (where \( 2^m \) is the larger power of 2 that divides \( 3x+1 \)).  The conjecture is: \emph{starting from any odd positive integer and repeating \( R(x) \) we eventually get to 1}. 

In \cite{Colussi} the set of odd positive integers \( x
 \) such that \( R^k(x) = 1 \) has been characterized as the set of odd integers whose binary representation belongs to a set of strings \( G_k \). Each string in \( G_k \) is the concatenation of \( k \) strings \( z_kz_{k-1}\dots z_1 \) where each \( z_i \) is a finite and contiguous extract from some power of a string \( s_i \) of length \( 2\cdot 3^{i-1} \) (the seed of order \( i \)). Clearly Collatz conjecture will be true if the binary representation of any odd integer belongs to some \( G_k \).
Lately Patrick Chisan Hew in \cite{Hew} showed that seeds \( s_i \) are the repetends of \( \frac{1}{3^i} \). 

Here two contributions to Collatz conjecture are given:
\begin{enumerate}
\item 
Collatz conjecture is expressed in terms of a function \( \rho(y) \) that operate on the set of all rational numbers \( 1/2 \leq y < 1 \) having finite binary representation. The main advantage of \( \rho(y) \) with respect to \( R(x) \) is that the denominator can be only \( 2 \) or \( 4 \) unlike \( R(x) \) whose demominator can be any power of 2.
\item 
We show that the binary representation of each odd positive integer \( x \) is a prefix of a power of infinitely many seeds \( s_i \) and we give an upper bound for the minimum \( i \) in terms of the length \( n \) of the binary representation of \( x \). 
\end{enumerate}
\end{abstract}

\maketitle

\section{Introduction}
The reduced Collatz function is defined as follows 
\begin{equation}\label{CollatzFunct}
R(x) = \frac{3x+1}{2^m}
\end{equation}
where \( x \) is an odd positive integer and \( 2^m \) is the highest power of 2 that divides \( 3x+1 \) (with no remainder). 

\FigureRho

Notice that \( m \geq 1 \) (since \( 3x+1 \) is even) and 3 does not divide \( R(x) \).
Collatz conjecture is
\begin{conjecture}\label{CollatzConj}
For all odd positive integer \( x \) there exists \( k \) such that \( R^{(k)}(x) = 1 \). 
\end{conjecture}
If \( R^{(k)}(x) = 1 \) for some \( k \) we say that \emph{\( x \) iterates to 1}. If \( R(x) \) iterates to 1 then \( x \) iterates to 1. Since \( R(x) \) is not a multiple of 3 in Conjecture \ref{CollatzConj} we can consider only the set \( X \) of odd positive integers \( x \) which are not a multiple of 3, i.e. such that \( x \bmod 3 \in \{1,2\} \). 

For \( x \in X \) let \( \psi(x) = x/2^n \) where \( n = \left\lceil \log_2 x \right\rceil \).  
If \( 1b_1\dots b_{n-2}1 \) is the binary representation of \( x \in X \) then \( 0.1b_1\dots b_{n-2}1 \) is the binary representation of \( \psi(x) \). By using the operator \( \sval{.} \) with the meaning "value of" we can write \( x = \sval{1b_1\dots b_{n-2}1} \) and 
\( \psi(x) = \sval{0.1b_1\dots b_{n-2}1} \).

Clearly \( \psi(x) \) is one-to-one and \( 1/2 \leq \psi(x) < 1 \). Let \( Y \) the set of all \( \psi(x) \) for \( x \in X \), i.e. the set of rational numbers \( 1/2 \leq y < 1 \) with finite binary representation \( y = \sval{0.1b_1\dots b_{n-2}1} \) and such that the odd positive integer \( y2^{n} = \sval{1b_1\dots b_{n-2}1} \) is not a multiple of 3, i.e. \( y \bmod (3\cdot 2^{-n}) \in \{2^{-n},2^{-n+1}\}\).

By using the one-to-one relation \( \psi(x) \) we can rewrite Collatz function as acting on \( y\in Y \)
\begin{equation}\label{CollatzFunctY}
\rho(y) = \psi\left(3\psi^{-1}(y)+1\right) 
= \frac{3y+2^{-n}}{2^{h}}
\end{equation}
where \( h = \left\lceil \log_2 (3y+2^{-n}) \right\rceil\). Notice that \( 1 < 3y+2^{-n} < 4 \) and so \( h = 1 \) or \( h = 2 \). 
More precisely \( h = 1 \) when \( 3y+2^{-n} < 2 \) (and so \( \rho(y) < 3/4 \)) and \( h = 2 \) when \( 3y+2^{-n} \geq 2 \) (and so \( \rho(y) \geq 3/4 \)). 
Of course \( \rho(y) \neq 3/4 \) otherwise \( y = 1 - 2^{-n}/3 \not\in Y \). 

Thus
\[ h = 1 \iff 3y+2^{-n} < 2 \iff \rho(y) > 3/4\] 
\[ h = 2 \iff 3y+2^{-n} \geq 2 \iff \rho(y) < 3/4\] 
Moreover \( \rho(y) = 1/2 \) iff \( 3y+2^{-n} = 2 \), i.e. iff \( y2^{n} = \frac{2^{n+1}-1}{3} \) and this can only happen 
for odd \( n \). For the particular case of \( n = 1 \) we have \( y = 1/2 \) 
and \( \rho(1/2) = 1/2 \).
Function \( \rho(y) \) is illustrated in
Figure \ref{FigureRho} both for even \( n \) and for odd \( n \).

The Collatz conjecture can be rewritten as
\begin{conjecture}\label{CollatzConjT}
For all \( y \in Y \) there exists \( k \) such that \( \rho^{(k)}(y) = 1/2 \). 
\end{conjecture}

\section{Some useful lemmas}\label{sect2} 
As usual we denote by \( q = \lfloor x/y \rfloor \) and \( r = x \bmod y \) the integer quotient \( q \) and the remainder  \( r \) of \( x/y \).

For irrational \( z > 0 \) and integer \( i \geq 0 \) let \( x_i = (iz) \bmod 1 \) the fractional part of \( iz \) and for \( k \geq 0 \) let \( X_k = \{ x_i : 0 \leq i < k\} \). 

Let represent  the elements of \( X_k \) in increasing order as points on a circle \( C \) of length \( 1 \) starting from point \( x_0 = 0 \) and moving counterclockwise (see  Figure \ref{XkCircle}). 
For \( x,x' \in X_k \) let \( d(x,x') \) the length of the counterclockwise arc from \( x \) to \( x' \). Notice that \( d(x,x') + d(x',x) = 1 \).

\begin{fact} \label{distance}
Let \( x_i,x_j \in X_k \). Then \( d( x_{i+t},x_{j+t}) = d( x_{i},x_{j}) \) for all \( t \) such that \( -\max(i,j) \leq t < k-\min(i,j) \). 
\end{fact}
\Proof Adding \( tz \) to \( x_i \) and \( x_j \) resorts to a rotation of length \( tz \) along the circle and so the distance does not change. The rotation is counterclockwise if \( t \geq 0 \) and clockwise if \( t \leq 0 \).
\qed
\medskip

\XkCircle

Points \( X_k \) split the circle into \( k \) arcs. 
Let \( \dmin{k} \) and \( \dmax{k} \) respectively the minimum and the maximum of the lengths of those arcs and let \( \ell \) the number of arcs of maximal length and \( s \) the number of arcs of minimal length. Moreover let \( q = \left\lfloor \dmax{k}/\dmin{k} \right\rfloor \) and \( r = \dmax{k} \bmod \dmin{k} \).

The case of \( k = 1 \) is uninteresting: there is only the point \( x_0 = 0 \) and \( \dmin{1} = \dmax{1} = 1 \). So from now on we assume \( k > 1 \).

\begin{fact} \label{minmaxonly}
Assume arcs of \( X_k \) take only the two lengths \( \dmin{k} < \dmax{k} \).  
 If \( k > 1 \) then \( r \neq 0 \). 
\end{fact}
\Proof 
If \( r = 0 \) then \( \dmax{k} = q \dmin{k} \) and \( \dmin{k} = 1/(\ell q + s) \) that contrasts \( z \) being irrational.
\qed
\medskip

\begin{lemma} \label{dmaxmin}
Assume arcs of \( X_k \) take only the two lengths \( \dmin{k} < \dmax{k} \). Also assume that one of the two arcs adjacent to \( x_0 \), say the arc \( (x_0,x_i) \),  has length \( \dmax{k} \) and the other, say the arc \( (x_0,x_j) \), has length \( \dmin{k} \). 

Point \( x_k \) splits \( (x_0,x_i) \) into one arc \( (x_k,x_i) \) of length \( \dmin{k} \) and one arc \( (x_0,x_k) \) of length \( \dmax{k} - \dmin{k} \). 

Moreover points \( x_{k+1},\dots, x_{k+\ell-1} \) split the same way the remaining \( \ell-1 \) arcs of length \( \dmax{k} \).
\end{lemma}

\Proof 
Point \( x_k \) can not split any arc \( (x_p,x_q) \) which is not adjacent to \( x_0 \) otherwise \( x_0 \) would split arc \( (x_{k-p},x_{k-q}) \) into two arcs both of length less than \( \dmax{k} \) (by Fact \ref{distance}) and that contrasts one arc adjacent to \( x_0 \) has length \( \dmax{k} \).

Point \( x_k \) can not split the arc \( (x_0,x_j) \) of length \( \dmin{k} \) otherwise \( (x_{k-j},x_0) \) would have length less than \( \dmin{k} \).

Thus \( x_k \) splits the arc \( (x_0,x_i) \) of length \( \dmax{k} \). Moreover the length of \( (x_k,x_i) \), which is the same as the length of \( (x_{k-i},x_0) \), should be \( \dmin{k} \) and the length of \( (x_0,x_k) \) is \( \dmax{k} - \dmin{k} \).

By Fact \ref{distance} the arcs of length \( \dmax{k} \) are exactly the arcs \( (x_t,x_{i+t}) \) for \( t = 0,\dots,k-i-1 \) for a total of \( \ell = k-i \) arcs.

Then, for \( t = 1,\dots, \ell-1 \), point \( x_{k+t} \) splits arc \( (x_t,x_{i+t}) \) into an arc \( (x_t,x_{k+t}) \) of the same length of arc \( (x_0,x_{k}) \) and an arc \( (x_{k+t},x_{i+t}) \) of the same length of the arc \( (x_{k},x_{i}) \).
\qed
\medskip

\begin{corollary} \label{archimede}
Assume that arcs of \( X_k \) take only two lengths \( \dmin{k} < \dmax{k} \) and the two arcs adjacent to \( x_0 \) have  different lengths.
Then:
\begin{enumerate}
\item for \( j = k + t\ell \) with \( 0 \leq t < q \),  arcs of \( X_j \) take only the two lengths \( \dmin{k} \) and \( \dmax{k}-t\dmin{k} \);

\item for \( k + t\ell < k' < k + (t+1)\ell \), arcs of \( X_{k'} \) take only the three lengths \( \dmin{k} \), \( \dmax{k}-t\dmin{k} \) and \( \dmax{k}-(t+1)\dmin{k} \);

\item for  \( k' = k + q\ell \), the set \( X_{k'} \) splits the circle into \( s' = \ell \) arcs of length \( \dmin{k'} = r \) and \( \ell' = \ell q + s  \) arcs of length \( \dmax{k'} = \dmin{k} \) and the two arcs adjacent to \( x_0 \) have  different lengths.
\end{enumerate}
\end{corollary}
\Proof
Use Lemma \ref{dmaxmin} \( q \) times.
 \qed
\medskip

Let define inductively \( k_h \) for all \( h \geq 0 \) by \( k_0 = 2 \) and \( k_h = k_{h-1} + q_{h-1}\ell_{h-1} \)
where \( q_{h-1} \) and \( \ell_{h-1} \) are the values of \( q \) and \( \ell \) for \( X_{k_{h-1}} \). Then \( k_h = 2+\sum_{j=0}^{h-1}q_{j}\ell_{j} \).

There are exactly two arcs of \( X_{k_{0}} = X_2 \) (of length \( \dmin{2} = \min(d,1-d) \) and \( \dmax{2} = \max(d,1-d) \) where \( d = z \bmod 1 \)). Moreover by Corollary \ref{archimede} 
arcs of \( X_{k} \) take only two lengths if \( k = k_{h} + t \ell_{h} \) for some \( h \) and \( 0 \leq t \leq q_{h}\) and take three different lengths otherwise. Moreover, \( \dmax{k} = \dmax{k_{h}}- t \dmin{k_{h}} \) for \( k_{h}+t \ell_{h} \leq k < k_{h}+(t+1) \ell_{h} \). 

\begin{lemma} \label{dmaxkhBound}
\( \frac{1}{k_{h}} < \dmax{k_{h}} < 
\frac{2}{k_{h}}
\).
\end{lemma}
\Proof Point 3 of Corollary \ref{archimede} implies
\(\ell_{h} = s_{h-1}+q_{h-1}s_{h} \) and so \( s_{h} = \frac{\ell_{h}-s_{h-1}}{q_{h-1}} < \frac{\ell_{h}}{q_{h-1}}\). 
Thus 
\[
\dmax{k_{h}} = \frac{1-s_{h} \dmin{k_{h}}}{\ell_{h}}< \frac{1}{\ell_{h}}  = \frac{\ell_{k}+s_{k}}{k_{h}\ell_{h}} < \frac{1+q_{h-1}}{k_{h}q_{h-1}} \leq \frac{2}{h}
\]
and
\[
\frac{1}{k_{h}} 
= \frac{\ell_{h}\dmax{k_{h}} + s_{h} \dmin{k_{h}}}{k_{h}} 
< \frac{(\ell_{h}+s_{h})\dmax{k_{h}}}{k_{h}} 
= \dmax{k_{h}}
\]
\qed
\medskip

The bounds from Lemma \ref{dmaxkhBound} holds only for those \( k \) such that \( k = k_h \) for some \( h \). The next three lemmas show what happen for \( k_{h} < k < k_{h+1} \).
 
\BoundDmax

\begin{lemma} \label{dmaxkBound} Let \( k = k_{h}+t \ell_{h} + j \) where \( 0 \leq t < q_{h} \) and \( 0 \leq j < \ell_{h} \). Then
\( \frac{1}{k_h}\cdot\frac{q_{h}-t}{q_{h}} < \dmax{k} < \frac{2}{k_{h}}\cdot \frac{q_{h}+1-t}{q_{h}+1} \).
\end{lemma}
\Proof
\[
\dmax{k} 
 = \dmax{k_h}-t\dmin{k_h}
 < \dmax{k_h}-\frac{t\dmax{k_{h}}}{q_{h}+1}
 < \frac{2}{k_h}\cdot\frac{q_{h}+1-t}{q_{h}+1}
\]
\[
\dmax{k} 
 = \dmax{k_h}-t\dmin{k_h}
 > \dmax{k_h}-\frac{t\dmax{k_{h}}}{q_{h}}
 > \frac{1}{k_h}\cdot\frac{q_{h}-t}{q_{h}}
\]
\qed
\medskip

Figure \ref{BoundDmax} show \( \dmax{k} \) in the range \( k_8 \leq k < k_9 \) for \( z = \log_2 3 \) compared to the bounds from Lemma \ref{dmaxkBound}. The range \( k_8 \leq k < k_9 \) is partitioned into \( q_{8} = 23 \) segments of \( \ell_{8} = 655 \) elements each. 

\begin{lemma} \label{epsilonkpkt} 
\( k\dmax{k} < 2\frac{q_{h}-t}{q_{h}}(1+t+\frac{j}{k_h}) 
\) for \( k_{h}\leq k < k_{h+1} \).
\end{lemma}
\Proof Let \( k = k_{h}+t \ell_{h} + j \) where \( 0 \leq t < q_{h} \) and \( 0 \leq j < \ell_{h} \).

By Point 2 of Corollary \ref{archimede}
\( \dmax{k} = \dmax{k_{h}}-t \dmin{k_{h}} \) and so 
\begin{eqnarray*}
k\dmax{k}
&=& (k_h+t\ell_{h}+j)(\dmax{k_{h}}-t \dmin{k_{h}})\\ 
&<& 
(k_h+t\ell_{h}+j)(1-\frac{t}{q_{h}})\dmax{k_{h}}
\\ 
&<&  
(k_h+tk_h+j)(1-\frac{t}{q_{h}})\frac{2}{k_h}
=
2\frac{q_{h}-t}{q_{h}}(1+t+\frac{j}{k_h})
\end{eqnarray*}
\qed
\BoundkDmax

\begin{lemma} \label{epsilonkpk} 
\( k\dmax{k} < \frac{q_{h}+5+5/q_{h}}{2} \) for \( k_{h}\leq k < k_{h+1} \). 
\end{lemma}
\Proof Let \( k = k_{h}+t \ell_{h} + j \) where \( 0 \leq t < q_{h} \) and \( 0 \leq j < \ell_{h} \).
The function \( f(t) = k\dmax{k}
= (k_h+t\ell_{h}+j)(\dmax{k_{h}}-t \dmin{k_{h}}) \) has a maximum for \( t_{max} = \frac{\dmax{k_{h}}
}{2\dmin{k_{h}}}-
\frac{(k_h+j)}{2\ell_{h}} \) and, taking account that \( \dmax{k_{h}} < 2/k_{h} \), \( \dmin{k_{h}} \leq \dmax{k_{h}}/q_{h} \) and \( \ell_{h} < k_{h} < 2\ell_{h} \), it follows that \( f(t_{max}) < \frac{q_{h}+5+5/q_{h}}{2} \).
\qed
\medskip

Figure \ref{BoundkDmax} shows \( k\dmax{k} \) in the range \( k_8 \leq k < k_9 \) for \( z = \log_2 3 \), the bound \( 2\frac{q_{h}-t}{q_{h}}(1+t+\frac{j}{k_h}) \) from Lemma \ref{epsilonkpkt} and the bound \( \frac{q_{h}+5+5/q_{h}}{2} \) from Lemma \ref{epsilonkpk}.%
\footnote{The bound \( \frac{q_{h}+5+5/q_{h}}{2} \) is \( \Theta(q_h) \). 
We can show that the \( q_j \) are the partial quotients of the continued fraction expansion \( z=q_0+\frac{1}{q_1+\frac{1}{q_2+\dots}} \).
 There are very few irrational \( z \) for which the quotients \( q_h \) are bounded by a constant. 
Unfortunately this seem not to hold for our case of interest \( z = \log_2 3 \).}

\section{Prefixes of seeds}
For each \( k \geq 0 \) let define the set
\[
P_k 
= 
\left\{ 
\frac{2^{\left\lceil i \log_2 3\right\rceil-1}}{3^i} 
: 0\leq i < k
\right\} 
\]
(Thinking in terms of binary representations we simply remove from \( 1/3^i \) the \( \left\lceil i \log_2 3\right\rceil-1 \) leading 0's after the dot.)
\begin{fact}
\( P_k \cap Y = \{1/2 \} \) for all \( k \geq 0 \). 
\end{fact}
\Proof 
Let \( p \in P_k \). If \( i = 0 \) then \( p = 1/2 \in Y \), otherwise \( p = \frac{2^{\left\lceil i \log_2 3\right\rceil-1}}{3^i} \) has infinite periodical binary representation and so \( p \not\in Y \).
\qed
\medskip

Take \( z = \log_2 3 \) as the irrational number used to build the set \( X_k \) in Section \ref{sect2} and let \( k_h = 2+\sum_{j=0}^{h-1}q_{j}\ell_{j} \). 

\begin{lemma} \label{denseP}
For all \( y \in Y \) and all \( \epsilon > 0 \)  there exists \( k \) and \( p \in P_{k} \) such that 
\( p \leq y < p+\epsilon \).
\end{lemma}
\Proof  The lemma is true for \( 1/2 \leq y < 1/2 + \epsilon \) since \( 1/2 \in P_{k} \) for all \( k \).

Let \( 1/2+\epsilon \leq y < 1 \); we are looking for \( i > 0 \) such that \( \frac{2^{\left\lceil i \log_2 3\right\rceil-1}}{3^i} \leq y < \frac{2^{\left\lceil i \log_2 3\right\rceil-1}}{3^i}+\epsilon \).

Let \( \delta = 1+\epsilon \) and let consider the strongest inequality \( \frac{2^{\left\lceil i \log_2 3\right\rceil-1}}{3^i} \leq y < \delta \frac{2^{\left\lceil i \log_2 3\right\rceil-1}}{3^i} \).
Taking the logarithms we obtain 
\[ 
\left\lceil i \log_2 3\right\rceil-1 - i\log_2 3  \leq \log_2 y  < \left\lceil i \log_2 3\right\rceil -1 - i\log_2 3 + \log_2 \delta 
\]
that can be rewritten as
\begin{equation}\label{inequality}
0 < -\log_2 y \leq (i\log_2 3) \bmod 1 < -\log_2 y + \log_2 \delta < 1
\end{equation}
(since \( \left\lceil i \log_2 3\right\rceil = 1 + \left\lfloor i \log_2 3\right\rfloor \), \( 1/2 \leq y < 1 \) and \( \delta \leq 1/2+y < 2y \).)

Let \( h \) minimum such that \( \frac{2}{k_h} < \log_2 \delta \) and take \( k=k_h \). Then \( d_{max}(k) < \log_2 \delta \) and there is at least one \( i < k \) that satisfy Inequality \ref{inequality}. Thus \( p \leq y < p+\epsilon \) for \( p = \frac{2^{\left\lceil i \log_2 3\right\rceil-1}}{3^i} \in P_{k} \).
\qed
\medskip



\begin{corollary} \label{distancePx}
Let \( 1/2= p_0 < p_1 < \dots < p_{k-1} < p_{k} = 1 \) the ordered sequence of \( p \in P_k \) with a last element \( p_k = 1 \) added. 
For all \( \epsilon > 0 \) there is \( k \) such that  \( p_{j+1}-p_j < \epsilon \) for all \( j = 0,\dots,k-1 \). 
\end{corollary}
\Proof Consequence of Lemma \ref{denseP}.
\qed
\medskip

\begin{lemma} \label{Prefixesx} The binary representation of each \( y \in Y \)  is a prefix of the binary representation of \( \frac{2^{\left\lceil i \log_2 3\right\rceil-1}}{3^i} \) for some \( i  \) (indeed for infinitely many \( i \)).
\end{lemma}
\Proof Let \( y = \sval{0.1b_1\dots b_{n-1}1} \). In Corollary 
\ref{distancePx} chose \( \epsilon =  2^{-n-1} \) and let \( j \) such that \( p_j \leq y < p_{j+1} \). Then \( y \) is a prefix of \( p_j  \) and \( p_j = \frac{2^{\left\lceil i \log_2 3\right\rceil-1}}{3^i} \) for some \( i < k \). 

Moreover there are infinitely many \( i \geq k \) such that  \( p_j \leq p= \frac{2^{\left\lceil i \log_2 3\right\rceil-1}}{3^i} < p_{j+1} \) and \( y \) is a prefix of all those \( p \). 
\qed
\medskip

\begin{theorem} \label{IntegerPrefixes} The binary representation of each odd integer \( x \)  is a prefix of some power of a seed \( s_i \) (whit leading 0's removed).
\end{theorem}
\Proof
 Let \( y = \psi(x) = x/2^n \) and let \( i \) such that the binary representation of \( y \) is a prefix of the binary representation of \( \frac{2^{\left\lceil i \log_2 3\right\rceil}}{3^i} \). Then \( x = \left\lfloor \frac{2^{\left\lfloor i \log_2 3\right\rceil+n}}{3^i} \right\rfloor\) is a prefix of a power of the seed \( s_i \).
\qed
\medskip

Lemma \ref{dmaxkhBound} and Lemma \ref{epsilonkpk} both give upper bounds for the minimum \( i \) such that the binary representation of \( x \) is a prefix of a power of \( s_i \).

The bound from Lemma \ref{dmaxkhBound} is \( k_h \) such that \( 2/k_h \leq 2^{-n-1} \) where \( n \) is the length of the binary representation of \( x \).

The bound from Lemma \ref{epsilonkpk} is \( k
\) such that \( k_h \leq k < k_{h+1} \) and \( \frac{q_h+5+5/q_h}{2k} \leq 2^{-n-1} \).


\section*{Bibliography}
\begin{thebibliography}{10}

\bibitem{Colussi}
Livio Colussi.
\newblock The convergence classes of Collatz function. 
\newblock Theor. Comp. Sci. 412 (2011)5409-5419.

\bibitem{Hew}
Patrik Chisan Hew.
\newblock Working in binary protects the repetends of \(1/3^h\): Comment on Colussi's 'The convergent classes of Collatz function'. 
\newblock Theor. Comp. Sci. 618 (2016)135-141.

\bibitem{Lagarias}
J.C. Lagarias.
\newblock The 3n+1 Problem: An Annotated Bibliography, II (2000-2009).
\newblock In http://arxiv.org/abs/math/0608208v5.


\end{thebibliography}
\end{document}